\documentclass[12pt]{amsart}

\usepackage{mathrsfs}
\usepackage{amsfonts}
\usepackage{amssymb}
\usepackage{amsfonts, amscd, amsmath, mathrsfs, amssymb, amsthm, amsxtra, bbding, epsfig, graphicx, latexsym, url, mathbbol, bbold}
\usepackage[papersize={8.3in,11.8in},textwidth=16.8cm,textheight=23.2cm,centering]{geometry}
\usepackage{enumerate}
\usepackage{caption} 

\usepackage{graphicx}
\usepackage{tikz}
\usepackage{float}

\usepackage{xcolor}
\definecolor{cite}{rgb}{0.00,0.60,1.00}
\definecolor{url}{rgb}{1.00,0.10,0.80}
\definecolor{link}{rgb}{0.00,0.00,1.00}
\usepackage[colorlinks,linkcolor=link,urlcolor=url,citecolor=cite,pagebackref,breaklinks]{hyperref}

\hypersetup{
pdfstartpage=1,
pdfstartview=FitH}

\DeclareFontFamily{U}{mathx}{\hyphenchar\font45}
\DeclareFontShape{U}{mathx}{m}{n}{
      <5> <6> <7> <8> <9> <10>
      <10.95> <12> <14.4> <17.28> <20.74> <24.88>
      mathx10
      }{}
\DeclareSymbolFont{mathx}{U}{mathx}{m}{n}
\DeclareMathAccent{\widecheck}{\mathalpha}{mathx}{"71}

\numberwithin{equation}{section}

\allowdisplaybreaks

\newtheorem{theorem}{Theorem}[section]
\newtheorem{lemma}{Lemma}[section]

\newtheorem{proposition}{Proposition}[section]

\makeatletter
\newcounter{roem}
\renewcommand{\theroem}{\Roman{roem}}

\newcommand{\c@org@eq}{}
\let\c@org@eq\c@equation
\newcommand{\org@theeq}{}
\let\org@theeq\theequation

\newcommand{\setroem}{
\let\c@equation\c@roem
 \let\theequation\theroem}

\newcommand{\setarab}{
\let\c@equation\c@org@eq
\let\theequation\org@theeq}
\makeatother

\newtheorem*{claim*}{Claim}

\theoremstyle{theorem}
\newtheorem{remark}{\bf Remark}

\newcommand{\ud}{\mathrm{d}}
\newcommand{\ue}{\mathrm{e}}

\newcommand{\kl}{\mathrm{Kl}}

\DeclareMathOperator*{\rint}{\hspace{-2mm}\ThisStyle{\footnotesize\hstretch{1.2}{\rotatebox{20}{$\SavedStyle\!\int\nolimits\!\!$}}}~\hspace{-2mm}}

\DeclareMathOperator*{\riiint}{\hspace{-2mm}\ThisStyle{\footnotesize\hstretch{1.2}{\rotatebox{20}{$\SavedStyle\!\int\!\!\!\!\!\!\!$}\rotatebox{20}{$\SavedStyle\!\int\!\!\!\!\!\!\!$}\rotatebox{20}{$\SavedStyle\!\int\!\!$}}}~\hspace{-2mm}}

\DeclareMathOperator*{\riiiiint}{\hspace{-2mm}\ThisStyle{\footnotesize\hstretch{1.2}{\rotatebox{20}{$\SavedStyle\!\int\!\!\!\!\!\!\!$}\rotatebox{20}{$\SavedStyle\!\int\!\!\!\!\!\!\!$}\rotatebox{20}{$\SavedStyle\!\int\!\!\!\!\!\!\!$}\rotatebox{20}{$\SavedStyle\!\int\!\!\!\!\!\!\!$}\rotatebox{20}{$\SavedStyle\!\int\!\!$}}}~\hspace{-2mm}}
\usepackage{scalerel}

\DeclareMathOperator{\Mod}{mod}

\renewcommand{\bmod}[1]{\,(\Mod{ #1})}

\newcommand{\bm}{\mathbf{m}}

\newcommand{\bA}{\mathbf{A}}

\newcommand{\bF}{\mathbf{F}}

\newcommand{\bR}{\mathbf{R}}
\newcommand{\bZ}{\mathbf{Z}}

\newcommand{\cM}{\mathcal{M}}

\newcommand{\cR}{\mathcal{R}}

\newcommand{\cV}{\mathcal{V}}
\newcommand{\cW}{\mathcal{W}}

\newcommand{\fA}{\mathfrak{A}}

\newcommand{\red}[1]{{\color{red}#1}}

\def\leq{\leqslant}

\def\geq{\geqslant}

\usepackage{graphicx}
\usepackage{tikz}

\begin{document}

\title{On the Brun--Titchmarsh theorem. II}

\author{Ping Xi}

\address{School of Mathematics and Statistics, Xi'an Jiaotong University, Xi'an 710049, P. R. China}
\address{Institute of Pure Mathematics, Xi'an Jiaotong University, Xi'an 710049, P. R. China}

\email{pingxi.cn@gmail.com, ping.xi@xjtu.edu.cn}

\author{Junren Zheng}

\address{School of Mathematics and Statistics, Xi'an Jiaotong University, Xi'an 710049, P. R. China}

\address{D{\'e}partment  de Math{\'e}matiques et Statistique,
Universit{\'e} de Montr{\'e}al, CP 6128 succ Centre-Ville, Montr{\'e}al, QC  H3C 3J7, Canada}
\email{junrenzheng03@gmail.com, junren@stu.xjtu.edu.cn}



\begin{abstract} 
Denote by $\pi(x;q,a)$ the number of primes $p\leqslant x$ with $p\equiv a\bmod q.$
We prove new upper bounds for $\pi(x;q,a)$ when $q$ is a large prime very close to $\sqrt{x}$, improving upon the classical work of Iwaniec (1982). The proof reduces to bounding a quintilinear sum of Kloosterman sums. To this end,  we introduce a new shifting argument inspired by Vinogradov--Burgess--Karatsuba, which goes beyond the classical Fourier-analytic approach. This argument relies on  a deep algebro-geometric result of Kowalski--Michel--Sawin on sums of products of Kloosterman sums.
\end{abstract}

\maketitle

\begin{flushright}
\dedicatory{{\it \small To Professor Yong-Gao CHEN on the occasion of his retirement}}
\end{flushright}
\setcounter{tocdepth}{1}


\section{Introduction and main results}
\label{sec:results}

\subsection{Main results}
We continue our study of the Brun--Titchmarsh theorem, i.e., bounding from above the 
counting function
\begin{align*}
\pi(x; q, a):=|\{p\leqslant x:p\equiv a\bmod q\}|,
\end{align*}
where $q\in\bZ^+$, $(a,q)=1$ and $x\geqslant3$. 
It is widely believed that 
\begin{align*}
\pi(x ; q, a)=(1+o(1)) \frac{1}{\varphi(q)}\frac{x}{\log x}
\end{align*}
as $q,x\rightarrow+\infty$ in a very broad range, but only the case
$q\leqslant (\log x)^A$ is known thanks to the classical Siegel--Walfisz theorem.

Throughout this paper, we assume $q\sim x^\varpi$ with $\varpi\in[0,1[$, and seek to determine a positive constant $C(\varpi)$ as small as possible, such that
\begin{align}\label{eq:BT}
\max_{(a,q)=1}\pi(x;q,a)\leqslant (C(\varpi)+\varepsilon)\frac{1}{\varphi(q)}\frac{x}{\log x}
\end{align}
holds for any $\varepsilon>0$ and $x>x_0(\varepsilon).$ Starting from the pioneering work of Titchmarsh \cite{Ti30}, various choices of $C(\varpi)$ have been obtained when $\varpi$ is in different ranges.
In particular, we should mention that Iwaniec \cite{Iw82} proved very strong bounds for $\pi(x;q,a)$ when $\varpi\in[9/20,2/3[.$ One of the main ingredients in his proof is the flexibility of the well-factorable remainder terms in the linear sieve, which he himself invented (see \cite{Iw80}). In particular,
Iwaniec showed that
\begin{align}\label{eq:Iwaniecconstant}
C(\varpi)=\frac{8}{6-7\varpi}
\end{align}
is admissible in \eqref{eq:BT}, for all $\varpi\in[9/20,2/3[.$ 
In \cite{XZ24}, we are able to improve the above bound of Iwaniec for all $\varpi\in[9/20,1/2[$ for general $q$, and for all $\varpi\in]1/2,4/7[$ when $q$ is a large prime.
More precisely, in the former case, we start from multiplicative characters, and the proof is then based on certain weighted moments of character sums or Dirichlet $L$-functions, to which the spectral theory of automorphic forms applies. 
In the latter case, however, we start from additive characters, and arrive at certain averages of Kloosterman sums with prime moduli. To capture sign changes of Kloosterman sums, we use combinatorial arguments and transfer the problem to the sums of products of many Kloosterman sums. Unfortunately, all the above methods 
fail to improve \eqref{eq:Iwaniecconstant} when $\varpi=1/2.$ 

In this paper, we may fill this gap for prime $q$, and obtain an improvement over \eqref{eq:Iwaniecconstant} when $\varpi$ is exactly at the threshold $1/2.$

\begin{theorem}\label{thm:varpi=1/2}
Let $q$ be a large prime.
Then we may take
\begin{align*}
C(\varpi)=\frac{240}{184-217\varpi}
\end{align*}
in \eqref{eq:BT} for all $\varpi\in[\frac{1}{2},\frac{34}{67}[.$ In particular, the choice $C(\frac{1}{2})=\frac{480}{151}$ is admissible.
\end{theorem}

Note that $\frac{480}{151}\approx 3.1788$ while Iwaniec \cite{Iw82} showed that one can take
$C(\frac{1}{2})=\frac{16}{5}=3.2$.
Theorem \ref{thm:varpi=1/2} is in fact an immediate consequence of the following general theorem by taking $\nu=8$.

\begin{theorem}\label{thm:varpi>1/2,small}
Let $q$ be a large prime.
For each positive integer $\nu\geqslant5,$ we may take
\begin{align*}
C(\varpi)=\frac{8}{6-7\varpi+\delta_\nu(\varpi)}
\end{align*}
in \eqref{eq:BT}, for all $\varpi\in[\frac{1}{2},\frac{\nu(2\nu+1)}{4\nu^2+\nu+4}[,$
where
\begin{align*}
\delta_\nu(\varpi)=\frac{2\nu-(3\nu+4)\varpi}{\nu(2\nu-1)}>0.
\end{align*}
\end{theorem}

\subsection{Ingredients of the proof}
To prove Theorem \ref{thm:varpi>1/2,small}, Iwaniec's version of linear sieve with a well-factorable remainder term is applied. By standard Fourier analysis, one then (essentially) arrives at the quintilinear sum
\begin{align}\label{eq:quintilinear-Kloosterman}
\sum_{k\leqslant K}\mathop{\sum\sum}_{n_1,n_2\sim N}\beta_{n_1}\overline{\beta_{n_2}}\mathop{\sum\sum}_{h_1,h_2\leqslant H}S(a(h_1\overline{n_1}-h_2\overline{n_2}),k;q),
\end{align}
where $S(\cdot,\cdot;q)$ is the Kloosterman sum defined by
\begin{align}\label{eq:Kloosterman}
S(m,n;q)=\sideset{}{^*}\sum_{x\bmod q}\ue\Big(\frac{mx+n\overline{x}}{q}\Big)
\end{align}
for $q\geqslant1$ (which is not necessarily a prime) and $m,n\in\bZ.$
The celebrated Weil's bound asserts that
\begin{align*}
|S(m,n;q)|\leqslant q^{\frac{1}{2}}(m,n,q)^{\frac{1}{2}}\tau(q).
\end{align*}
This exhibits the square-root cancellations in \eqref{eq:Kloosterman}, which is best possible in the individual sense. Iwaniec employed the Weil's bound to control the above quintilinear sum of Kloosterman sums. To beat Iwaniec's barrier, we need to capture cancellations among averages of Kloosterman sums to arrive at stronger inequalities.
Note that, in Iwaniec's arguments, one is allowed to take $H\approx N=x^{\frac{1}{2}-\varepsilon}q^{-\frac{3}{4}}$ and $K=q^{2+\varepsilon}/x$. When $q$ is significantly larger than $x^{\frac{1}{2}}$, the above $K$ would have a reasonably large size, and we may combine the variables $h_1,h_2,n_1,n_2$ to create bilinear/trilinear sums of Kloosterman sums. This was carried out in our previous paper \cite{XZ24}. However, for $q\approx x^{\frac{1}{2}}$, we find $H\approx N\approx q^{\frac{1}{4}}$ and $K\approx q^{\varepsilon}$ in the above choices, and cancellations from the sum over $k$ are no longer possible if one follows the previous arguments.

Suppose that we are in the critical situation that 
\begin{align}\label{eq:critical}
q\sim x^{\frac{1}{2}},\ \ H=N=q^{\frac{1}{4}},\ \ K=q^{\varepsilon}.
\end{align}
For each $k\leqslant q^\varepsilon$, we now encounter quadrilinear sums of Kloosterman sums, each of the summations is of length $q^{\frac{1}{4}}$. 
To deal with such extremely short averages, we are inspired by a {\it shifting trick} of Vinogradov~\cite{Vi58}, Burgess~\cite{Bu62}  and Karatsuba~\cite{Ka68}, which has been  developed and explained more recently by Friedlander–Iwaniec~\cite{FI85}, Fouvry–Michel~\cite{FM98} and Kowalski–Michel–Sawin~\cite{KMS17, KMS20}.  We introduce a variant of this method that exploits the special additive structure in \eqref{eq:quintilinear-Kloosterman}, allowing us to go beyond the limitation of the direct Fourier-analytic methods. 
More precisely, we perform the shift
$(h_1,h_2)\mapsto(h_1+rs_1,\;h_2+rs_2)$ for $r,s_1,s_2\in\bZ.$ Averaging over the three new variables in suitable ranges, we introduce 
$u := ak(h_1\overline{n}_1 - h_2\overline{n}_2), 
v :=ak(s_1\overline{n}_1 - s_2\overline{n}_2)$, and apply H\"{o}lder's inequality to transfer the problem to a high moment of averages of Kloosterman sums (Lemma \ref{lm:Kloostermanmoment-2dim}). For this we appeal to a deep study of Kowalski--Michel--Sawin \cite{KMS20} on the stratification method for {\it sum-product} of Kloosterman sheaves.
Finally, we are able to detect cancellations among Kloosterman sums in \eqref{eq:quintilinear-Kloosterman} in the 
critical situation \eqref{eq:critical}.
See also \cite{FSX25} for our recent applications to bounding trilinear and quadrilinear character sums.

\subsection*{Notation and conventions}  
We write $\ue(t)=\ue^{2\pi\mathrm{i}t}$. We use $\varepsilon$ to denote a small positive number, which might be different at each occurrence. The notation $n\sim N$ means $N<n\leqslant2N.$ Denote by $\varphi$ the Euler function, and by $\tau_k$ the $k$-fold divisor function. 
For a function $g\in L^1(\bR),$ define the Fourier transform by
\begin{align*}
\widehat{g}(y)=\rint_{\bR} g(x)\ue(-yx)\ud x.
\end{align*}

\subsection*{Acknowledgements}
We thank \'Etienne Fouvry, Philippe Michel and the two referees for valuable suggestions.
This work is supported in part by NSFC (No. 12025106) and Shaanxi NSF (No. 2025JC-QYCX-002).

\smallskip

\section{From sieves to exponential sums}
\label{sec:sievemethods}

Let $y=x(\log x)^{-3}$ and introduce a smooth function $\varPhi$ which is supported on $[y,x+y]$, satisfying
\begin{align*}
\begin{cases}
\varPhi(t)=1\ \ &\text{for } t\in[2y,x],\\
\varPhi(t)\geqslant0\ \ &\text{for } t\in\bR,\\
\varPhi^{(j)}(t)\ll_j y^{-j} &\text{for all } j\geqslant0,
\end{cases}
\end{align*}
and $\widehat{\varPhi}(0)=\int_\bR \varPhi(t)\ud t=x.$
For $(d, q)=1$ put
\begin{align}\label{eq:A(d)}
A(d)
:= \sum_{\substack{n\equiv a\bmod q\\d\mid n}} \varPhi(n)=
\sum_{n\equiv a\overline{d}\bmod q} \varPhi(dn)
\end{align}
and
\begin{align}\label{eq:r(d)}
r(d)
&:=A(d)-\frac{x}{dq}.
\end{align}
For $z>2$, denote by $P_q(z)$ the product of primes less than $z$ and coprime to $q$.

We now recall the following sieve inequality in \cite[Proposition 4.1]{XZ24}.
\begin{proposition}\label{prop:sieveinequality}
Let $\varepsilon>0$ sufficiently small and $D<x^{1-\varepsilon}.$ For $z=\sqrt{x/2}$ and any $M,N\geqslant1$ with $MN=D,$ we have
\begin{align}\label{eq:pi(x;q,a)-sievebound} 
\pi(x; q, a)
\leqslant \{2+O(\varepsilon)\}\frac{x}{\varphi(q)\log D}
+\sum_{t\leqslant T(\varepsilon)}\mathop{\sum_{m\leqslant M}\sum_{n\leqslant N}}_{mn\mid P_q(z)}\alpha_t(m)\beta_t(n) r(mn),
\end{align}
where $T(\varepsilon)$ depends only on $\varepsilon,$ $|\alpha_t(m)|,|\beta_t(n)|\leqslant 1,$ and $r(\cdot)$ is defined by $\eqref{eq:A(d)}$ and $\eqref{eq:r(d)}.$
\end{proposition}

From Poisson summation it follows that
\begin{align*}
A(d)
=\frac{1}{dq}\sum_{h\in\bZ} \widehat{\varPhi}\Big(\frac{h}{dq}\Big)\ue\Big(\frac{ah\overline{d} }{q}\Big).
\end{align*}
From the rapid decay of $\widehat{\varPhi}$, we may write, for $d\leqslant D$, that
\begin{align*}
r(d)
&=A(d)-\frac{x}{dq}
=\frac{1}{dq} \sum_{0<|h|\leqslant H}\widehat{\varPhi}\Big(\frac{h}{dq}\Big)
\mathrm{e}\Big(\frac{ah\overline{d}}{q}\Big)+O(x^{-10})
\end{align*}
with $H:=x^{\varepsilon-1}MNq$.

From Proposition \ref{prop:sieveinequality}, together with a dyadic device, it suffices to prove that
\begin{align}\label{eq:R(M,N)-expectedbound}
\cR(M,N):=\sum_{m\sim M}\sum_{n\sim N} \frac{\alpha_m\beta_n}{mn} \sum_{h\leqslant H}\widehat{\varPhi}\Big(\frac{h}{mnq}\Big)
\mathrm{e}\Big(\frac{ah\overline{mn}}{q}\Big)\ll \frac{x}{(\log x)^{2025}},
\end{align}
where $|\alpha_m|,|\beta_n|\leqslant 1,$ $M,N\geqslant1$ with $MN=D$ chosen as large as possible. Note that the restriction $(mn,q)=1$ is always kept in mind henceforth.
We transform the above multiple exponential sum in $\cR(M,N)$ into averages of complete Kloosterman sums.

By Cauchy--Schwarz inequality, we obtain
\begin{align*} 
\cR(M,N)^2
&\leqslant \frac{1}{M}\sum_{m\in\bZ}W\Big(\frac{m}{M}\Big)\Big|\sum_{h\leqslant H}\sum_{n\sim N} \frac{\beta_n}{n}\widehat{\varPhi}\Big(\frac{h}{mnq}\Big)
\ue\Big(\frac{ah\overline{mn}}{q}\Big)\Big|^2\\
&=\frac{1}{M}\mathop{\sum\sum\sum\sum}_{h_1,h_2\leqslant H,~n_1,n_2\sim N}\frac{\beta_{n_1}\overline{\beta_{n_2}}}{n_1n_2}\sum_{m\in\bZ}W\Big(\frac{m}{M}\Big)\widehat{\varPhi}\Big(\frac{h_1}{mn_1q}\Big)
\overline{\widehat{\varPhi}\Big(\frac{h_2}{mn_2q}\Big)}
\ue\Big(\frac{a(h_1\overline{n_1}-h_2\overline{n_2})\overline{m}}{q}\Big),
\end{align*}
where $W$ is a fixed smooth function which dominates the indicator function of $[1,2].$
By Poisson summation, the $m$-sum becomes
\begin{align*} 
\frac{1}{q}\sum_{k\in\bZ}\widetilde{\varPhi}(k;h_1,h_2,n_1,n_2,q)S(a(h_1\overline{n_1}-h_2\overline{n_2}),k;q),
\end{align*}
where $S(\cdot,\cdot;q)$ is the Kloosterman sum defined by \eqref{eq:Kloosterman}, and
\begin{align}\label{eq:varPhi-tilde}
\widetilde{\varPhi}(k;h_1,h_2,n_1,n_2,q)=\rint_{\bR}W\Big(\frac{t}{M}\Big)\widehat{\varPhi}\Big(\frac{h_1}{tn_1q}\Big)
\overline{\widehat{\varPhi}\Big(\frac{h_2}{tn_2q}\Big)}\ue\Big(\frac{-kt}{q}\Big)\ud t.
\end{align}
From integration by parts, we find
\begin{align*} 
\widetilde{\varPhi}(k;h_1,h_2,n_1,n_2,q)\ll x^2M(1+|k|M/q)^{-A}
\end{align*}
for any $A\geqslant0.$ Hence we may truncate the $k$-sum at $K:=M^{-1}q^{1+\varepsilon}$ with a negligible error term. Moreover, we bound the
contributions from $k=0$ by
\begin{align*} 
&\ll \frac{x^2}{qN^2}\mathop{\sum\sum\sum\sum}_{h_1,h_2\leqslant H,~n_1,n_2\sim N}(h_1n_2-h_2n_1,q)\\
&\ll \frac{x^2}{qN^2}\mathop{\sum\sum}_{l_1,l_2\leqslant 2HN}\tau_2(l_1)\tau_2(l_2)(l_1-l_2,q)\\
&\ll x^{2+\varepsilon}HN^{-1}(q^{-1}HN+1).
\end{align*}
For the remainder of the paper, we assume that $q$ is a large prime and that $M>2q^{\varepsilon}$, which implies $K<q/2$. Therefore,
\begin{align}\label{eq:cR(M,N)^2-varSigma}
\cR(M,N)^2
&\leqslant \varSigma+O(x^{2+\varepsilon}HN^{-1}(q^{-1}HN+1))
\end{align}
with
\begin{align} \label{eq:varSigma}
\varSigma
&=\frac{1}{\sqrt{q}M}\mathop{\sum\sum\sum\sum}_{h_1,h_2\leqslant H,~n_1,n_2\sim N}\frac{\beta_{n_1}\overline{\beta_{n_2}}}{n_1n_2}\sum_{1\leqslant |k|\leqslant K}\widetilde{\varPhi}(k;h_1,h_2,n_1,n_2,q)\kl(ak(h_1\overline{n_1}-h_2\overline{n_2}),q),
\end{align}
where we write
\begin{align}\label{eq:Kloosterman-kl}
\kl(z,q)=\frac{S(z,1;q)}{\sqrt{q}}.
\end{align}

\smallskip

\section{Proof of Theorem \ref{thm:varpi>1/2,small}}

We are now ready to prove Theorem \ref{thm:varpi>1/2,small}. To bound $\varSigma$ (as defined by \eqref{eq:varSigma}) from above, we first introduce an auxiliary function
\begin{align*} 
\phi(\lambda)=
\begin{cases}
\min\{\lambda,1,H+ 1-\lambda\},& \text{for~}\lambda\in[0,H],\\
0,&\text{otherwise}.
\end{cases}
\end{align*} 
This allows us to write
\begin{align*} 
\varSigma
&=\frac{1}{\sqrt{q}M}\mathop{\sum\sum\sum\sum}_{h_1,h_2\in\bZ,~n_1,n_2\sim N}\frac{\phi(h_1)\phi(h_2)\beta_{n_1}\overline{\beta_{n_2}}}{n_1n_2}\sum_{1\leqslant |k|\leqslant K}\widetilde{\varPhi}(k;h_1,h_2,n_1,n_2,q)\kl(ak(h_1\overline{n_1}-h_2\overline{n_2}),q).
\end{align*}
Note that integration by parts implies
\begin{equation}\label{eq:phi-fourier}
\widehat{\phi}(\xi)\ll \min\{H,|\xi|^{-1},|\xi|^{-2}\}.
\end{equation}

\subsection{The shifting trick}
For all non-zero integers $r,s_1,s_2$, we find
\begin{align*} 
\varSigma
&=\frac{1}{\sqrt{q}M}\mathop{\sum\sum}_{n_1,n_2\sim N}\frac{\beta_{n_1}\overline{\beta_{n_2}}}{n_1n_2}\sum_{1\leqslant |k|\leqslant K}\mathop{\sum\sum}_{h_1,h_2\in\bZ}\phi(h_1+rs_1)\phi(h_2+rs_2)\\
&\ \ \ \ \ \times \widetilde{\varPhi}(k;h_1+rs_1,h_2+rs_2,n_1,n_2,q)\kl(ak((h_1+rs_1)\overline{n_1}-(h_2+rs_2)\overline{n_2}),q).
\end{align*}
In what follows, we would like to take advantage of summing over $r,s_1,s_2$ non-trivially. To this end, suitable separations of variables are necessary.
Note that
\begin{align*} 
\widetilde{\varPhi}&(k;h_1+rs_1,h_2+rs_2,n_1,n_2,q)
=\rint_{\bR}W\Big(\frac{t}{M}\Big)\widehat{\varPhi}\Big(\frac{h_1+rs_1}{tn_1q}\Big)
\overline{\widehat{\varPhi}\Big(\frac{h_2+rs_2}{tn_2q}\Big)}\ue\Big(\frac{-kt}{q}\Big)\ud t\\
&=\riiint_{\bR^3}W\Big(\frac{t}{M}\Big)\varPhi(\xi_1)\varPhi(\xi_2)\ue\Big(-\frac{(h_1+rs_1)\xi_1}{tn_1q}+\frac{(h_2+rs_2)\xi_2}{tn_2q}-\frac{kt}{q}\Big)\ud t\ud \xi_1\ud\xi_2.
\end{align*}
Making the changes of variables $\xi_1\mapsto \xi_1n_1/s_1$ and $\xi_2\mapsto \xi_2n_2/s_2$, $\widetilde{\varPhi}(\cdots)$ becomes
\begin{align*} 
&\frac{n_1n_2}{s_1s_2}\riiint_{\bR^3}W\Big(\frac{t}{M}\Big)\varPhi\Big(\frac{\xi_1n_1}{s_1}\Big)\varPhi\Big(\frac{\xi_2n_2}{s_2}\Big)\ue\Big(-\frac{(h_1/s_1+r)\xi_1}{tq}+\frac{(h_2/s_2+r)\xi_2}{tq}-\frac{kt}{q}\Big)\ud t\ud \xi_1\ud\xi_2.
\end{align*}
Also, for $j=1,2,$
\begin{align*} 
\phi(h_j+rs_j)
&=\rint_{\bR}\widehat{\phi}(\eta) \ue((h_j+rs_j)\eta)\ud\eta=\rint_{\bR}\widehat{\phi}\Big(\frac{\eta}{s_j}\Big) \ue((h_j/s_j+r)\eta)\frac{\ud\eta}{s_j}
\end{align*}
by Fourier inversion. It is also useful to write
\begin{align*} 
(h_1+rs_1)\overline{n_1}-(h_2+rs_2)\overline{n_2}=
h_1\overline{n_1}-h_2\overline{n_2}+(s_1\overline{n_1}-s_2\overline{n_2})r.
\end{align*}
We now sum over $r\sim R,$ $s_1,s_2\sim S$ with 
\begin{align}
R,S\geqslant 1,\ \ 4RS=H,
\end{align}
and we find 
\begin{align*} 
\varSigma
&\ll\frac{1}{\sqrt{q}MRS^6}\mathop{\sum\sum\sum\sum\sum\sum}_{|h_1|,|h_2|\leqslant 2H,~n_1,n_2\sim N,~s_1,s_2\sim S}\sum_{1\leqslant |k|\leqslant K}\riiiiint_{\bR^5}W\Big(\frac{t}{M}\Big)\Big|\widehat{\phi}\Big(\frac{\eta_1}{s_1}\Big)\widehat{\phi}\Big(\frac{\eta_2}{s_2}\Big)\Big|\varPhi\Big(\frac{\xi_1n_1}{s_1}\Big)\varPhi\Big(\frac{\xi_2n_2}{s_2}\Big)\\
&\ \ \ \ \ \times \Big|\sum_{r\sim R}\theta_{r;t,\eta_1,\eta_2,\xi_1,\xi_2}\kl(ak(h_1\overline{n_1}-h_2\overline{n_2})+ak(s_1\overline{n_1}-s_2\overline{n_2})r,q)\Big|\ud t\ud\eta_1\ud\eta_2\ud\xi_1\ud\xi_2
\end{align*}
with $\theta_{r;t,\eta_1,\eta_2,\xi_1,\xi_2}=\ue((\xi_2-\xi_1)r/(tq)+(\eta_1+\eta_2)r).$
From \eqref{eq:phi-fourier} it follows that
\begin{align*} 
\varSigma
&\ll\frac{1}{\sqrt{q}MRS^6}\mathop{\sum\sum\sum\sum\sum\sum}_{|h_1|,|h_2|\leqslant 2H,~n_1,n_2\sim N,~s_1,s_2\sim S}\sum_{1\leqslant |k|\leqslant K}\\
&\ \ \ \ \ \times\riiiiint_{\bR^5}W\Big(\frac{t}{M}\Big)\min\Big\{H,\frac{S}{|\eta_1|},\frac{S^2}{\eta_1^2}\Big\}\min\Big\{H,\frac{S}{|\eta_2|},\frac{S^2}{\eta_2^2}\Big\}\varPhi\Big(\frac{\xi_1n_1}{s_1}\Big)\varPhi\Big(\frac{\xi_2n_2}{s_2}\Big)\\
&\ \ \ \ \ \times \Big|\sum_{r\sim R}\theta_{r;t,\eta_1,\eta_2,\xi_1,\xi_2}\kl(ak(h_1\overline{n_1}-h_2\overline{n_2})+ak(s_1\overline{n_1}-s_2\overline{n_2})r,q)\Big|\ud t\ud\eta_1\ud\eta_2\ud\xi_1\ud\xi_2
\end{align*}
Recall the definition of $W$ and $\varPhi,$ we infer
\begin{align*} 
\varSigma
&\ll\frac{x^{2+\varepsilon}}{\sqrt{q}N^2RS^2}\mathop{\sum\sum\sum\sum\sum\sum}_{|h_1|,|h_2|\leqslant 2H,~n_1,n_2\sim N,~s_1,s_2\sim S}\sum_{1\leqslant |k|\leqslant K}\\
&\ \ \ \ \ \times \Big|\sum_{r\sim R}\theta_r \kl(ak(h_1\overline{n_1}-h_2\overline{n_2})+ak(s_1\overline{n_1}-s_2\overline{n_2})r,q)\Big|
\end{align*}
for some coefficient $\theta_r$ with $|\theta_r|=1.$ Moreover, we may impose additional restrictions $h_1n_2\neq h_2n_1$ and $s_1n_2\not\equiv s_2n_1\bmod q$, with an error $O(x^{2+\varepsilon}q^{-\frac{3}{2}}(NS)^{-1}H^2K(NS+q))$.

\subsection{Applying Hölder's inequality}
Put
\begin{align*}
\rho(b,c)=\mathop{\sum\sum\sum\sum\sum\sum\sum}_{\substack{|h_1|,|h_2|\leqslant 2H,~n_1,n_2\sim N,~s_1,s_2\sim S,~1\leqslant |k|\leqslant K\\
ak(h_1\overline{n_1}-h_2\overline{n_2})\equiv b\bmod{q}\\ ak(s_1\overline{n_1}-s_2\overline{n_2})\equiv c\bmod{q}\\ h_1n_2\neq h_2n_1}}1
\end{align*}
for $b,c\bmod q$ with $(c,q)=1$.
Hence
\begin{align*}
\varSigma
&\ll\frac{x^{2+\varepsilon}}{\sqrt{q}N^2RS^2}\mathop{\sum\sum}_{\substack{b,c\bmod q\\ (c,q)=1}}\rho(b,c)\Big|\sum_{r\sim R}\theta_r
\kl(b+rc,q)\Big|+x^{2+\varepsilon}q^{-\frac{3}{2}}(NS)^{-1}H^2K(NS+q).
\end{align*}
By H\"older's inequality, for each $\nu\in\bZ^+,$
\begin{align}\label{eq:varSigma-Holder}
\varSigma
&\ll\frac{x^{2+\varepsilon}}{\sqrt{q}N^2RS^2}\varSigma_1^{1-\frac{1}{\nu}}(\varSigma_2\varSigma_3)^{\frac{1}{2\nu}}
\end{align}
with
\begin{align*}
\varSigma_j=\mathop{\sum\sum}_{\substack{b,c\bmod q\\ (c,q)=1}}\rho(b,c)^j,\ \ j=1,2,
\end{align*}
and
\begin{align*}
\varSigma_3=\mathop{\sum\sum}_{\substack{b,c\bmod q\\ (c,q)=1}}\Big|\sum_{r\sim R}\theta_r
\kl(b+rc,q)\Big|^{2\nu}.
\end{align*}

Trivially we have
\begin{align*}
\varSigma_1\ll H^2N^2S^2K.
\end{align*}

To bound $\varSigma_2,$ we appeal to the following lemma.

\begin{lemma}\label{lm:varSigma2-upperbound}
With the above notation, we have
\begin{align*}
\varSigma_2
&\ll \Big(1+\frac{HN^3K}{q}\Big)\Big(1+\frac{SN^3K}{q}\Big)\Big(1+\frac{H}{N}\Big)\Big(1+\frac{S}{N}\Big)(HNS)^2Kq^{\varepsilon}.
\end{align*}
\end{lemma}

\proof Write $T_1=32HN^3K/q$ and $T_2=32SN^3K/q.$
Hence 
\begin{align}\label{eq:varSigma2-varSigma2(E,E')}
\varSigma_2\ll \log^2N\sup_{1\leqslant E,E'\leqslant N}\varSigma_2(E,E'),
\end{align}
where
$\varSigma_2(E,E')$ denotes the number of tuples
$$(e,e',h_1,h_2,h_1',h_2',n_1,n_2,n_1',n_2',s_1,s_2,s_1',s_2',k,k',t_1,t_2)\in \bZ^{18}$$
with
\begin{gather}
e\sim E,~e'\sim E',~0\leqslant|h_1|,|h_2|,|h_1'|,|h_2'|\leqslant 2H,~n_1,n_2\sim N/E,~n_1',n_2'\sim N/E',\nonumber
\\
s_1,s_2,s_1',s_2'\sim S,~1\leqslant |k|,|k'|\leqslant K,~0\leqslant |t_1|\leqslant T_1/(EE'),~0\leqslant |t_2|\leqslant T_2/(EE'),\nonumber
\\
ke'n_1'n_2'(h_1n_2-h_2n_1)=k'en_1n_2(h_1'n_2'-h_2'n_1')+t_1q,\label{eq:varSigma2-t1}\\ 
ke'n_1'n_2'(s_1n_2-s_2n_1)=k'en_1n_2(s_1'n_2'-s_2'n_1')+t_2q,\label{eq:varSigma2-t2}\\
(h_1n_2-h_2n_1)(h_1'n_2'-h_2'n_1')(s_1n_2-s_2n_1)(s_1'n_2'-s_2'n_1')\neq0,\nonumber\\
(n_1,n_2)=(n_1',n_2')=1.\nonumber
\end{gather}

We consider the contributions from $t_1t_2\neq0$, and those from $t_1t_2=0$ can be treated easily by divisibility. We look at \eqref{eq:varSigma2-t1} firstly. Given an integer $\ell$ with $0<|\ell|\leqslant4HNK/E,$ consider the equations
\begin{align}\label{eq:h1h2n1n2''''t-equation1}
k(h_1 n_2-h_2 n_1)=\ell,
\end{align}
and
\begin{align}\label{eq:h1h2n1n2''''t-equation2}
\ell e'n_1'n_2'=k'en_1n_2(h_1'n_2'-h_2'n_1')+t_1q.
\end{align}
Note that 
\eqref{eq:h1h2n1n2''''t-equation2} is equivalent to the congruence equation
\begin{align}\label{eq:h1h2n1n2''''t-equation3}
k'en_1n_2(h_1'n_2'-h_2'n_1')\equiv -t_1q\bmod{e'n_1'n_2'}.
\end{align}
Given $e',h_1',h_2',n_1',n_2',t_1$, the number of tuples $(k',e,n_1,n_2)$ satisfying \eqref{eq:h1h2n1n2''''t-equation3}
with $1\leqslant |k'|\leqslant K,$ $e\sim E, n_1,n_2\sim N/E$ is 
\begin{align*}
\ll \Big(\frac{KN^2}{E}\frac{(h_1'n_2'-h_2' n_1',d'n_1'n_2')}{e'n_1'n_2'}+1\Big)N^{\varepsilon}
\ll \Big(\frac{E'K}{E}(h_1'n_2'-h_2' n_1',d'n_1'n_2')+1\Big)N^{\varepsilon},\end{align*}
where the factor $N^\varepsilon$ comes from the bound $\tau_4(m)\ll m^\varepsilon.$
After determining $k',e,n_1,n_2$ in \eqref{eq:h1h2n1n2''''t-equation3}, the variable $\ell$ is also determined as shown by \eqref{eq:h1h2n1n2''''t-equation2}.
We now turn to \eqref{eq:h1h2n1n2''''t-equation1}, from which we see that 
the number of pairs $(h_1,h_2,k)$ with $0\leqslant|h_1|,|h_2|\leqslant 2H$ and $1\leqslant |k|\leqslant K$ is $\ll N^\varepsilon(1+H/n_1)\ll N^\varepsilon(1+EH/N).$ 
With such choices of $e,e',n_1,n_2,n_1',n_2',k,k',$ we find the number of tuples $(s_1,s_2,s_1',s_2',t_2)$ satisfying \eqref{eq:varSigma2-t2} is at most
\begin{align*}
\ll \frac{S^2}{EE'}(1+T_2)\Big(1+\frac{S}{n_1}\Big)\ll \frac{S^2}{EE'}(1+T_2)\Big(1+\frac{ES}{N}\Big).
\end{align*}

Combining all above arguments, we arrive at
\begin{align*}
\varSigma_2(E,E')
&\ll q^{\varepsilon} (1+T_1)(1+T_2)\Big(\frac{S}{EE'}\Big)^2\Big(1+\frac{EH}{N}\Big)\Big(1+\frac{ES}{N}\Big)\\
&\ \ \ \ \times \mathop{\sum\sum\sum\sum\sum}_{\substack{e'\sim E',~|h_1'|,|h_2'|\leqslant 2H,~n_1',n_2'\sim N/E'\\ h_1'n_2'\neq h_2'n_1'}} \Big(\frac{E'K}{E}(h_1'n_2'-h_2'n_1',e'n_1'n_2')+1\Big)\\
&\ll q^{\varepsilon}(1+T_1)(1+T_2)\Big(\frac{S}{EE'}\Big)^2\Big(1+\frac{EH}{N}\Big)\Big(1+\frac{ES}{N}\Big)\Big(\frac{1}{E}+\frac{1}{E'}\Big)(HN)^2K.
\end{align*}
Now the lemma follows by recalling \eqref{eq:varSigma2-varSigma2(E,E')}.
\endproof

\begin{remark}
To illustrate the strength of Lemma $\ref{lm:varSigma2-upperbound},$ it is helpful to remark that the bounds for $\Sigma_1$ and $\Sigma_2$ roughly coincide in the critical range $q\sim x^{\frac{1}{2}},\ S \leq H \approx N \approx
 x^{\frac{1}{4}},\ K = x^{o(1)}$. 
\end{remark}

\subsection{Sums of products of Kloosterman sums}
To estimate $\varSigma_3$, we appeal to the following deep result due to Kowalski, Michel and Sawin \cite{KMS20}. 

\begin{lemma}\label{lm:Kloostermanmoment-2dim}
Let $q$ be a large prime and $M\geqslant3.$
For each positive integer $\nu,$ and any sequence of complex coefficients $(\alpha_m)$ with $|\alpha_m|\leqslant1,$ we have
\begin{align*}
\mathop{\sum\sum}_{\substack{b,c\bmod q\\ (c,q)=1}}\Big|\sum_{m\leqslant M}\alpha_m
\kl(b+mc,q)\Big|^{2\nu}
\ll q^2 M^\nu+q M^{2\nu}.
\end{align*}
\end{lemma}

\proof
Denote by $\fA$ the sum in question and write $\cM=[1,M]^{2\nu}.$
Opening the power and by a change of variable, we get
\begin{align*}
|\fA|\leqslant \sum_{\bm\in\cM}|\fA(\bm)|,
\end{align*}
where for $\bm=(m_1,m_2,\cdots,m_{2\nu})\in\cM,$
\begin{align*}
\fA(\bm)=\mathop{\sum\sum}_{\substack{b,c\bmod q\\ (c,q)=1}}\prod_{1\leqslant j\leqslant 2\nu}\kl((b+m_j)c,q).
\end{align*}
According to \cite[Theorem 4.6]{KMS20}, there exist affine varieties
\begin{align*}
\cV\subseteq \cW\subseteq \bA_{\bZ}^{2\nu}
\end{align*}
defined over $\bZ$ such that
\begin{align*}
\mathrm{codim}(\cV)=\nu, \quad \mathrm{codim}(\cW) \geqslant \nu/2,
\end{align*}
and
\begin{align*}
\fA(\bm)\ll
\begin{cases}
q^2, \ \ & \text { if } \bm\in \cV(\bF_q),\\
q^{\frac{3}{2}},& \text { if } \bm\in (\cW-\cV)(\bF_q),\\
q, & \text { if } \bm\not\in \cW(\bF_q).
\end{cases}
\end{align*}
According to Kowalski, Michel and Sawin \cite[Theorem 4.5]{KMS20}, one may refer to a counting lemma of Xu \cite[Lemma 1.7]{Xu20}, getting
\begin{align*}
|\cV(\bF_q)\cap\cM|\ll M^\nu,\quad |(\cW-\cV)(\bF_q)\cap\cM|\ll M^{\frac{3\nu}{2}}.
\end{align*}
Therefore,
\begin{align*}
|\fA|
&\leqslant \sum_{\bm\in\cV(\bF_q)\cap\cM}|\fA(\bm)|+\sum_{\bm\in(\cW-\cV)(\bF_q)\cap\cM}|\fA(\bm)|+\sum_{\bm\in\cM\setminus\cW(\bF_q)}|\fA(\bm)|\\
&\ll q^2 M^\nu+q^{\frac{3}{2}} M^{\frac{3\nu}{2}}+q M^{2\nu}\\
&\ll q^2 M^\nu+q M^{2\nu}.
\end{align*}
This established the lemma.
\endproof

By Lemma \ref{lm:Kloostermanmoment-2dim}, we also find
\begin{align*}
\varSigma_3\ll q^2 R^\nu+q R^{2\nu}.
\end{align*}
With the above bounds for $\varSigma_1,\varSigma_2$ and $\varSigma_3,$ it follows from \eqref{eq:varSigma-Holder} that
\begin{align*}
\varSigma
&\ll\frac{x^{2+\varepsilon}}{\sqrt{q}N^2RS^2}(H^2N^2S^2K)^{1-\frac{1}{2\nu}}\Big(1+\frac{HN^3K}{q}\Big)^{\frac{1}{2\nu}}\Big(1+\frac{SN^3K}{q}\Big)^{\frac{1}{2\nu}}\\
&\ \ \ \times\Big(1+\frac{H}{N}\Big)^{\frac{1}{2\nu}}\Big(1+\frac{S}{N}\Big)^{\frac{1}{2\nu}}(q^2 R^\nu+ q R^{2\nu})^{\frac{1}{2\nu}}+x^{2+\varepsilon}q^{-\frac{3}{2}}(NS)^{-1}H^2K(NS+q).
\end{align*}

\subsection{Concluding Theorem \ref{thm:varpi>1/2,small}}
Suppose that 
\begin{align}\label{eq:Hlarge}
H>4q^{\frac{1}{\nu}}.
\end{align}
Taking 
\begin{align*}
R=q^{\frac{1}{\nu}},\ \ S=\tfrac{1}{4}Hq^{-\frac{1}{\nu}},
\end{align*}
we have
\begin{align*}
\varSigma
&\ll x^{2+\varepsilon}q^{-\frac{1}{2}+\frac{1}{2\nu}+\frac{1}{\nu^2}}H^{2-\frac{2}{\nu}}N^{-\frac{1}{\nu}}K^{1-\frac{1}{2\nu}}\Big(1+\frac{HN^3K}{q}\Big)^{\frac{1}{2\nu}}\Big(1+\frac{SN^3K}{q}\Big)^{\frac{1}{2\nu}}\Big(1+\frac{H}{N}\Big)^{\frac{1}{\nu}}\\
&\ \ \ \ +x^{2+\varepsilon}q^{-\frac{3}{2}}(NS)^{-1}H^2K(NS+q).
\end{align*}
It now follows from \eqref{eq:cR(M,N)^2-varSigma} that
\begin{align*}
\cR(M,N)^2
&\ll x^{2+\varepsilon}H^{2-\frac{2}{\nu}}N^{-\frac{1}{\nu}}K^{1-\frac{1}{2\nu}}q^{-\frac{1}{2}+\frac{1}{2\nu}+\frac{1}{\nu^2}}\Big(1+\frac{HN^3K}{q}\Big)^{\frac{1}{2\nu}}\Big(1+\frac{HN^3K}{q^{1+\frac{1}{\nu}}}\Big)^{\frac{1}{2\nu}}\Big(1+\frac{H}{N}\Big)^{\frac{1}{\nu}}\\
&\ \ \ \ +x^{2+\varepsilon}q^{-1}HN^{-1}(HN+q)(1+q^{-\frac{1}{2}+\frac{1}{\nu}}K).
\end{align*}
Recalling that $H=x^{\varepsilon-1}MNq$ and $K=M^{-1}q^{1+\varepsilon}$, we arrive at
\begin{align*}
\cR(M,N)^2
&\ll x^{\frac{2}{\nu}+\varepsilon}M^{1-\frac{3}{2\nu}}N^{2-\frac{3}{\nu}}q^{\frac{5}{2}-\frac{2}{\nu}+\frac{1}{\nu^2}}\Big(1+\frac{N^4q}{x}\Big)^{\frac{1}{2\nu}}\Big(1+\frac{N^4q^{1-\frac{1}{\nu}}}{x}\Big)^{\frac{1}{2\nu}}\Big(1+\frac{Mq}{x}\Big)^{\frac{1}{\nu}}\\
&\ \ \ \ +x^{1+\varepsilon}qM(1+x^{-1}MN^2)(1+M^{-1}q^{\frac{1}{2}+\frac{1}{\nu}}).
\end{align*}

To guarantee \eqref{eq:R(M,N)-expectedbound}, and in view of \eqref{eq:Hlarge}, it suffices to require that both
\begin{gather*}
x^{\frac{2}{\nu}} M^{1-\frac{3}{2\nu}} N^{2-\frac{3}{\nu}} q^{\frac{5}{2}-\frac{2}{\nu}+\frac{1}{\nu^2}}
\Big(1+\frac{N^{4} q}{x}\Big)^{\frac{1}{2\nu}}
\Big(1+\frac{N^{4} q^{1-\frac{1}{\nu}}}{x}\Big)^{\frac{1}{2\nu}}
\Big(1+\frac{M q}{x}\Big)^{\frac{1}{\nu}}
< x^{2-\varepsilon'}\\
x\,q\,M\Big(1+\frac{MN^{2}}{x}\Big)\Big(1+\frac{q^{\frac{1}{2}+\frac{1}{\nu}}}{M}\Big)
< x^{2-\varepsilon'}
\end{gather*}
for some $\varepsilon'=\varepsilon'(\varepsilon)>0$.

To guarantee \eqref{eq:R(M,N)-expectedbound} and in view of \eqref{eq:Hlarge}, we require that
\begin{gather*}
x^{\frac{2}{\nu}}M^{1-\frac{3}{2\nu}}N^{2-\frac{3}{\nu}}q^{\frac{5}{2}-\frac{2}{\nu}+\frac{1}{\nu^2}}\Big(1+\frac{N^4q}{x}\Big)^{\frac{1}{2\nu}}\Big(1+\frac{N^4q^{1-\frac{1}{\nu}}}{x}\Big)^{\frac{1}{2\nu}}\Big(1+\frac{Mq}{x}\Big)^{\frac{1}{\nu}}<x^{2-\varepsilon'}\red{,}\\
xqM(1+x^{-1}MN^2)(1+M^{-1}q^{\frac{1}{2}+\frac{1}{\nu}})<x^{2-\varepsilon'}
\end{gather*}
with a suitably small $\varepsilon'$ in terms of $\varepsilon.$
The two simultaneous constraints are satisfied if we take
\begin{align*}
M=x^{1-\varepsilon'}q^{-1},\ \ N=\min\{x^{\frac{2\nu-1}{2(2\nu-3)}}q^{-\frac{3\nu^2-\nu+2}{2\nu(2\nu-3)}},~x^{\frac{\nu}{2\nu-1}}q^{-\frac{3\nu^2+2}{2\nu(2\nu-1)}}\}x^{-\varepsilon'}.
\end{align*}
In particular, for $\frac{1}{2}\leqslant\varpi<\tfrac{\nu(2\nu+1)}{4\nu^2+\nu+4},$
we may simplify the choice of $N$ to
\begin{align*}
N=x^{\frac{\nu}{2\nu-1}-\varepsilon'}q^{-\frac{3\nu^2+2}{2\nu(2\nu-1)}},
\end{align*}
in which case we obtain
\begin{align*}
D=MN=x^{\frac{3\nu-1}{2\nu-1}-2\varepsilon'}q^{-\frac{7\nu^2-2\nu+2}{2\nu(2\nu-1)}}.
\end{align*}
We also need to take $ \nu\geq 5 $ to guarantee \eqref{eq:Hlarge}.
Substituting the above level $D$ in Proposition \ref{prop:sieveinequality}, we conclude Theorem \ref{thm:varpi>1/2,small}.

\smallskip

\bibliographystyle{plain}

\end{document}